\documentstyle[12pt]{article}
\newcommand{\Title}[1]{{\Large \bf \begin{center} {#1}
  \end{center}}\vspace*{1mm}}
\newcommand{\Author}[1]{{\bf \begin{center} {#1}
  \end{center}}\vspace*{1mm}}
\newtheorem{theorem}{Theorem}
\newtheorem{lemma}{Lemma}          
\newtheorem{corollary}{Corollary}
\begin{document}

\Title{Graphs as rotations}

\Author{Dainis ZEPS
\footnote{This research is supported by Institution of Mathematics and
Computer Science of Latvian University under grant 93.602and by KAM of Caroline University by grant 
GAUK 351, GA\v{C}R 2167.} \footnote{Author's address: Institute of Mathematics and Computer Science, 
University of Latvia, 29 Rainis blvd., Riga, Latvia.
{dainize@mii.lu.lv}}
}
\begin{abstract}
Using a notation of corner between edges when graph has a  fixed rotation, 
i.e. cyclical order of edges around vertices, we define combinatorial 
objects - combinatorial maps as pairs of permutations, 
one for vertices and one for faces. Further, we define multiplication 
of these objects, 
that coincides with the multiplication of permutations. We consider 
closed under multiplication classes of combinatorial maps that 
consist of closed classes of combinatorial maps with fixed edges 
where each such class is defined by a knot. One class among
 them is special, containing selfconjugate maps.
\end{abstract} 

\section{Introduction}

Usually,  speaking about graphs on 
orientable surfaces \cite{whi 73}, we  suppose that edges incident with a vertex are  
cyclically ordered, or, in other words,  the rotation of such graphs is fixed.
 Actually the objects that are fixed are ordered pairs 
of edges that follow one another in the 
cyclical order around a vertex. We name such pair of edges a {\em corner} due to a 
natural geometrical  interpretation. Now thinking in terms of corners 
 we may say, that just the corners 
around vertices are ordered cyclically. Moreover, 
the cyclical sequences of corners of different vertices do not overlap, 
so they form together one common permutation on the set of  all 
corners.  We look on this permutation as the rotation of 
this graph. Further, if we "forget" the other 
information about the graph except this one 
permutation, then the graph cannot be 
restored uniquely. Only if we take permutations both of the 
graph and its dual graph (on this surface), then the graph is restorable (up to 
isomorphism) from them, but, of course,  together with its rotation, i. e. its imbedding in the surface. 
Generally, two arbitrary permutations $p_1$ and $p_2$, 
whenever all cycles of $p_1^{-1} \cdot p_2$ are 
transpositions, define some graph with fixed 
rotation(up to isomorphism). In this work we are dealing with combinatorial 
objects called {\em (combinatorial) maps}, which are pairs of permutations.

Combinatorial maps have been studied in many works by several authors, for 
example \cite{vin 83,sta 80,sta 83,sta 88,tut 84,lit 88,BonLit 95}. This work 
is done independently from them.

Many ground things considered in this work are  
well known for many years. But in some points we differ from the attitudes of the named authors. 
For example, we define a binary operation on  the maps 
that corresponds to multiplication of permutations. This 
operation is defined only if the product $p_1^{-1} 
\cdot p_2$ is equal for both the graphs and so it is common for 
the set of graph that is closed under this operation. 
This set is as big as the set of all permutations (of some 
fixed degree). Consequently it gives possibility to establish 
one-one map between all maps and all permutations. We 
are arrived at the situation where one permutation 
corresponds to one graph, and, if we could consider 
this permutation as the rotation of this graph, then this 
graph is given, when just this rotation is given.

This insight gives us right to speak now about graphs as 
rotations in the place of graphs with rotation, when we 
consider them as a whole set, that is closed under one 
common operation.

All permutations of fixed order comprise the class of maps closed under multiplication of maps.

In chapter \ref{knot} we enter a notion of combinatorial knot that is uniquely connected with 
the edges of that map. It corresponds to zigzag walks in \cite{BonLit 95}. 
In \cite{Liu 94} it is called a closed travel obtained by the T.T.Rule
. Every knot defines some subclass of combinatorial maps these with  fixed edges. 
Besides, we define selfconjugate map with as many links in its knot as edges in it.  
We prove that the subclass of selfconjugate maps is a subgroup of the fixed class of combinatorial 
maps and that every map can be expressed as the multiplication of the knot of a map 
with some selfconjugate map. 

\section{Permutations}
{\em Fixed elements} in permutations are these that form cycles of length one. 
{\em Transpositions} are cycles of length two. 
Permutations with only fixed elements and 
transpositions we  call {\em involutions}. 
Involutions without fixed elements we call {\em matchings}.

Our permutations act on a universal set $C$ the elements of which we call sometimes corners (or 
simply elements). For a permutation $P$ and $c \in C$  $c^P$ denotes that element of C to which
$c$ goes over. We are multiplying permutations from left to right.  For two 
permutations $P$ and $Q$  $P^Q$ equals to $Q^{-1} 
\cdot P \cdot Q$ and is the conjugate permutation to $P$ (with respect to $Q$). When $S$ 
is a cycle $(c_1 c_2 \ldots )$ in $P$, then $S^Q$ is the cycle $(c_1^Q c_2^Q \ldots )$ in 
$P^Q$. (For arbitrary permutations $P, Q, R$ it holds: $(P^Q)^{R}=P^{Q \cdot R}$ and $(P 
\cdot Q)^{R} = P^R \cdot Q^R$.) $I$ denotes identical permutation, i.e. that with all elements fixed.
 One way to express a cycle 
$(c_1 c_2 \ldots c_k) (k > 0)$ of some permutation as 
the multiplication of transpositions is this: $(c_k c_{k-1}) \ldots (c_2 
c_1)$.

\section{Combinatorial maps}

We say that two permutations $P$ and $Q$ of equal order $n$ are
{\em differing} whenever $i^P \ne i^Q$ for every $i \in [1 \ldots n]$. 
A pair of permutations $(P, Q)$ 
we call {\em combinatorial map} whenever $P$ and $Q $ are differing 
permutations and  $P^{-1} \cdot Q$ 
is a matching.
   
We denote the product $P^{-1} \cdot Q$ with the letter $\pi$. It is 
convenient to choose one fixed $\pi$.
In  that case we are speaking about a map $ P$, keeping $Q$ and matching $\pi$ in mind, or
saying, that $\pi$ is fixed. The definition of the combinatorial map is such, 
that we do not need to speak 
about the orders of the permutations. But sometimes we do, and then it is easy 
to see, that the order of the 
permutations $P, Q, \pi$ is even, let us put it equal to $n = 2 \times m$.

 We say, that $P, Q,  \pi$ act on
$2m$ {\em corners} or elements. If $C$ is the  set of corners , then it is
divided by $\pi$  into $m$ pairs. Let $(c_1, c_2)$ be such  a pair, i.e.
$c_1^{\pi} = c_2$. Then there exists such pair of corners $(c, c')$ that
$(c,c') \cdot P = P \cdot (c_1, c_2)$. We call the pair $(c, c')$ an 
{\em edge of the map} and $(c_1, c_2)$ we call a {\em next edge of the
map}.  We call matching $\pi$ the {\em next-edge-matching(n-matching) of a map}. For a map
$(P , P \cdot \pi)$ there exists some other unique matching $\varrho$ satisfying
$\varrho \cdot P = P \cdot \pi$. We call matching $\varrho$ {\em
edge-matching(e-matching) of the map}. We say that $\pi$ {\em contains} next edges and
$\varrho$ contains edges of the map $P$ , n-matching $\pi$ being fixed.

For a combinatorial map $(P, Q)$ we call the map $(Q, P)$ its {\em dual} combinatorial 
map. We use also denotation $\bar P$ for $Q$, saying that $\bar P$ is dual to $P$ by $\pi$ being fixed.

\section{Closed classes of combinatorial maps}

Let us define the
multiplication of two maps $S$ and $T$ (with fixed the same n-matching
$\pi$) as usual multiplication of permutations, i.e. we put 
$$ (S_1 , S_2) \cdot (T_1, T_2) = (S_1 \cdot T_1, S_1 \cdot T_1 \cdot \pi)$$ 
by definition, where on the left side of the expression sign $'\cdot '$
stands for multiplication of maps.

It is easy to see that 
$$(S_1 , S_2) \cdot (T_1, T_2) = (S_1 \cdot  S_2, S_1 \cdot T_2)$$
In practice, writing $S \cdot T$ we mean both
multiplication of permutations $S$ and $T$ and multiplication of
corresponding maps with common n-matching $\pi$. 

It is easy to see that it hold 
$$ 1) \ \ \ (S , {\bar S}  ) \cdot (T, {\bar T}) = (S \cdot T , S \cdot  {\bar T}) $$  
$$ 2) \ \ \ (S , {\bar S}  ) \cdot ({\bar T}, T) = (S \cdot {\bar T} , S \cdot  T) $$
$$ 3) \ \ \ ({\bar S} , S  ) \cdot (T, {\bar T}) = ({\bar S} \cdot T , {\bar S} \cdot T) $$
$$ 4) \ \ \ ({\bar S} , S  ) \cdot ({\bar T}, T) = ({\bar S} \cdot {\bar T} , {\bar S} \cdot T) $$

Multiplication for
maps with different n-matchings is not defined. 

 For any fixed n-matching (with m next edges) the set of all permutations $P_{2m}$ 
with $2m$ elements defines and forms
a class of maps that is closed against multiplication of maps. 
For two different permutations we have correspondingly two different (not equal) maps.
Hence, we have established one-one 
map between a set of permutations and a set of maps.
So, we may think in terms of one fixed n-matching and a set of 
all permutations, under each permutation seeing a graph $P$ 
that stands for the ordered pair $(P, {\bar P})$. So, for arbitrary 
permutation $P$ we have also a graph $(P, {\bar P})$).
 Let us take the map $(I, \pi)$ as identity map. Graphically interpreted it
consists of $m$ isolated edges. Its dual map
 $(\pi,I)$ consists of m isolated loops.

$P^{-1}$ being reverse permutation
of $P$, $( P^{-1}, P^{-1} \cdot \pi)$ is called {\em reversed} map of $(P, {\bar P})$. So, both
map and its reverse map belong to the same class of maps closed under
multiplication. Hence, because all permutations with $2m$ elements comprise  the 
symmetric group $S_{2m}$  \cite{and 72, wie 64}, so similarly,
 all maps with one fixed n-matching form a
group that is isomorphic to $S_{2m}$. 

If $(a, b)$ is a transposition, $(a, b) \cdot P$ can 
be graphically interpreted as a union 
of two vertices at the given corners $a$ and $b$ (or split of the  vertex which contains both corners)in 
the graph corresponding to the map $P$. Remembering that 
each permutation can be expressed as a sequence of multiplications of simple
transpositions, each map can be generated by this sequence of the
operations of union of two corners (or split of vertex): 
if $P = (a_l b_l)(a_{l-1} b_{l-1})...(a_1 b_1)$ 
for some $l > 0$
then $P_0 = I$, $P_k = (a_k,b_k) \cdot P_{k-1}$, $0 < k \leq l$, $P = P_l$.
Graphically, this shows l unions of vertices (or splits)at given corners giving the graph $P$.

 It is convenient to  choose one special class of
maps with n-matching equal to $\pi = (1 2) (3 4) \ldots \\  (2i-1 \ \ 2i) \ldots 
(2m-1 \ \ 2m)$. The maps of this class we call {\em normal}. It is
convenient for practice to get used to some particular chosen
n-matchings and this one named  normal is sufficiently natural.

The multiplication of maps gives right to the following graphical 
interpretation. The right member of multiplication is the map is the 
given one, i.e. that is subjected to changes, but the left comprise 
these transpositions that must affect the right one. 
Consequently, a map can be considered as move from one map to some 
other. In category language a map is both the object and the arrow. 

In applications  it is convenient to choose one big natural number, say $M$,  and consider it as the bound 
for the possible number of edges in  maps . Taking this number M as the order 
of all permutations, defining maps, then formally all 
maps have $M$ edges, but most of them should be isolated. Such 
view is convenient also because it can be used directly in computer application. 
The only pay for 
this approach is that we loose possibility to have isolated edges, when we 
actually need them, from 
some graph theoretical point of view. In the graphical interpretation then 
we have, say, at 
the beginning $M$ isolated edges (being initialized), and each permutation, multiplied 
to the currant graph, is some set of 
transpositions, that union and split vertices at distinct corners, giving in this way  new graphs.

\section{Classes of maps with fixed edges} 

For a map $P$ with
n-matching $\pi$ its e-matching $\varrho$ is equal to $P \cdot \pi
\cdot P^{-1} (=\pi^{P^{-1}})$. Except $P$ there exist other maps with the same e-matching. Truly, every
permutation $Q$ against which matchings $\pi$ and $\varrho$ are
conjugate, defines a map with the same e-matching $\varrho$. 

Let us define a set
$K_{\varrho}$ of maps as class of all maps whose e-matching is
$\varrho$. For different $\varrho$ these classes $K_{\varrho}$ are subclasses of maps of all
class of maps with fixed $\pi$, under which one this subclass is
special, namely, $K_{\pi}$, that comprise the maps with e-matching
equal to $\pi$. This subclass is not empty, because the maps $(\pi, I)$
and $(I , \pi)$ belong to it. Really, $\varrho_I = \pi^I = \pi$.

\begin{theorem}\label{fix1}
 For two maps $S$, $T$ with n-matching $\pi$ e-matching of
their multiplication $S \cdot T$ is equal to $ \varrho_T^{S^{-1}}$,
i.e. $\varrho_{S \cdot T}^S = \varrho_T$.
\end{theorem}
{\bf Proof} 
$$\varrho_{S \cdot T}
 = S \cdot T \cdot \pi \cdot T^{-1}
\cdot S^{-1} =
S \cdot \varrho_T \cdot S^{-1} = \varrho_T^{S^{-1}}.$$

   Let us denote by $P \cdot K_{\sigma}$ class of maps $\{P \cdot Q |
Q \in K_{\sigma}\}$. From theorem \ref{fix1} we know that $P \cdot K_{\sigma}$  
goes into $K_{\varrho}$, where $\varrho = \sigma^{P^{-1}}$.  We are going to
prove the equality of these classes. First we start with this. 

\begin{theorem}
 $K _{\pi}$ is a group.
\end{theorem}
{\bf Proof}  If $\varrho_S = \varrho_T
= \pi$ then $\varrho_{S \cdot T} = S \cdot T \cdot \pi \cdot T^{-1}
\cdot S^{-1} =
S \cdot \pi \cdot S^{-1} = \pi$.  Further, $\varrho_{S^{-1}} = S^{-1}
\cdot  \pi \cdot S = \pi$. Identity map also belongs to this
subclass , and, so, it forms a subgroup. 

\begin{theorem}\label{fix3}
 $P \cdot K_{\sigma} = K_{\sigma ^{P^{-1}}}$.  $(P \cdot K_{T^P} = K_T.)$
\end{theorem}
{\bf Proof} If $\sigma = \pi $ then $P \cdot K_{\sigma}$
is a left coset of $K_{\sigma}$ equal to $K_{\pi^{P^{-1}}}$ and theorem is
right.

  Let $\sigma \ne \pi$ and $Q$ be such map that $Q \cdot K_{\pi}
= K_{\sigma}$, i.e. $K-{\sigma}$ is a left coset of $K_{\pi}$
 and $Q$ is some of its elements: $\varrho_Q = \pi^{Q^{-1}}$ and $Q$ belongs to $K_{\sigma}$ by
theorem~1.  If $P \cdot Q = R$ then we have $P \cdot K_{\sigma} = P
\cdot Q \cdot K_{\pi} = R \cdot K_{\pi} = K_{\pi^{R^{-1}}} = K_{\sigma^{P^{-1}}}$.
As a conclusion of this theorem we have that the class $K_{\pi}$ and
its left cosets are classes with fixed edges. 

\begin{corollary} 
Left cosets of $K_{\pi}$ are classes of maps with fixed edges.
\end{corollary}

 So arbitrary map $P$ with
e-matching $\varrho = \pi^{P^{-1}}$ belongs to $K_{\varrho} = P \cdot K
_{\pi}$.  

For two involutions $\sigma$ and $\tau$ ($\sigma^2 = \tau^2 =
I$) we write $\sigma \subseteq \tau$, if every transposition of $\sigma$ is
also a transposition in $\tau$.  

\begin{lemma}
 $P^{\pi} = P $ iff $P \in
K_{\pi}$.
\end{lemma}
{\bf Proof} $\varrho_P = \pi \equiv \pi = P \cdot \pi \cdot
P^{-1} \equiv P = \pi \cdot P \cdot \pi \equiv P = P^{\pi}$. 

 Let $c$
be a cycle of $P$. Cycle $c^{\pi}$ is called its {\em conjugate }cycle (with
respect to $\pi$).  If $c = c^{\pi}$ then it is called a {\em
selfconjugate} cycle. If every cycle in $P$ has its conjugate cycle in
$P$ or is selfconjugate then $P$ is called {\em selfconjugate map}.
Lemma says that $K_{\pi} $ is the class of selfconjugate maps with respect to $\pi$.

   We reformulate this in theorem. 

\begin{theorem}
 $K_{\pi} $ is the class of
selfconjugate maps with respect to $\pi$.
\end{theorem}

  We can reveal the structure of these maps as follows.  

\begin{theorem}
 $K_{\pi}$ is isomorphic to $S_m \cdot S_2^m$.
\end{theorem}
{\bf Proof}   Let $P \in K_{\pi}$ and $P$ work on $C$ and $C_1 \cup C_2$ be
one partition of $C$ induced by $\pi$. Then there is such involution
$\sigma \subseteq \pi$ that $P = Q \cdot \sigma$ and cycles of $Q$
belong completely to $C_1 $ or $C_2$. i.e. if a cycle  $c$ goes into $C_1$ then
$c^{\pi}$ goes into $C_2$ or reversely. $Q$ can be expressed as $Q_1
\cdot Q_2$ where $Q_1$ has corners of $C_1$ and $Q_2$ corners of $C_2$. 
Then $Q_1$ and $Q_2$ are both isomorphic  between themselves and to some 
permutation in $S_m$ and $\sigma$ isomorphic to  some permutation of 
$S_2^m$.

\begin{theorem}
 $ \mid K_{\pi} \mid  = m !  \times  2^m$.
\end{theorem}
{\bf Proof}   $ \mid S_m \mid = m!  $ ;    $ \mid S_2^m \mid  = 2^m$. 

\begin{theorem}
 $K_{\pi}$ has (together with itself) $(2 m - 1)!!$ left cosets.
\end{theorem}  
{\bf Proof}    We have as many left cosets of $K_{\pi}$ as many
e-matchings can be generated. Their number is $(2 m -1)!!$ .

One more way to convince oneself that all left cosets of $K_{\pi}$ have 
the same number of elements, i.e. $(2m-1)!!$, and that 
they really do not overlay each other, is to see that  the equality 
$(2m-1)!! \times m! 2^m = (2m)!$  holds.

\section{Combinatorial knot}\label{knot}

Let $(P , P \cdot \pi)$ be a map on set of corners $C$ and $\varrho = \pi^{P^{-1}}$.
Let $C_1 \cup C_2$ be some partition of $C$ such that it is induced both
by $\pi$ and $\varrho$, i.e. for every both edge and next edge one of
its corner belongs to $C_1$ and other to $C_2$. In this case we say
that $C$ is {\em well } partitioned or {\em well colored in two colors}
or we say that $C_1 \cup C_2$ is {\em well coloring} of $C$ induced by
this map $(P , P \cdot \pi)$. 

\begin{theorem}\label{kno1}
 There always exists well coloring of
$C$ induced by arbitrary map $(P , P \cdot \pi)$.
\end{theorem}  
{\bf Proof} \ \ Let $c_1 c_2 \ldots
c_{2k}$ be a cycle of corners such that $c_2 = c_1^{\pi}$ and $c_3 =
c_2^{\varrho}$ and so on in alternating way, i.e. $c_{2i} =
c_{2i-1}^{\pi}$ and $c_{2i+1} = c_{2i}^{\varrho}$,for $i = 1, \ldots
k$, where $c_{2k}^{\varrho} = c_{2k+1} = c_1$. 
Let us suppose for an instant that $c_{2k-1}^{\pi} = c_{2k} = c_1$.
Then $c_{2k-1} = c_1^{\pi} = c_2$ and $c_{2k-2}^{\varrho} =c_2$ 
and $c_{2k-2} = c_2^{\varrho} = c_3$ and so on,
 until $c_{k+1} = c_k$, but it isn't possible.
It follows, that these cycles may have only even number of elements.

Then we may put odd
elements of cycle in $C_1$ and even elements in $C_2$ . If this
cycle runs through all corners, then we have only one possibility to
color all corners. Otherwise we choose arbitrary non colored corner taking
it as $c_1$ and proceed as before. In the end we get well partition
$C_1 \cup C_2$ induced by the map $(P , P \cdot \pi)$.  

Let us define
permutation $\mu$ having cycles as in the previous proof, i.e. if $C_1
\cup C_2$ is well coloring of the set $C$, then  $c^{\mu} = c^{\pi}$ 
if $c \in  C_1$ and $c^{\mu} = c^{\varrho}$ if $c \in C_2$. Permutation 
$\mu$ is called {\em combinatorial knot}. In place we use shorter name {\em knot}
 for  $\mu$. In a graphical 'corner' interpretation $\mu$ really is 
(alternating) knot  of the graph \cite{BurZie 85}. In \cite{BonLit 95} it is called zigzag walk.

\begin{lemma}
 If $\mu$ is a knot then if $\mu'$ is obtained with some
cycle of $\mu$ changed in the opposite  direction then $\mu'$ is also a
knot.
\end{lemma}
{\bf Proof} \ \ 
 Starting a new cycle in the proof of \ref{kno1} we can choose a
corner arbitrary . Consequently, cycle in a knot can go in one or another direction. 
It follows that if $\mu $ is a knot then also $\mu 
^{-1}$ is a knot.  It is easy to see that a knot $\mu$ depends only on $\pi$
and $\varrho$, so it is common for all class $K_{\varrho}$. The following theorem shows that, if we
consider a map $(\pi, \mu)$ then $\varrho_{\mu} = \varrho$. 

\begin{theorem}
$\pi^{\mu} = \varrho$.
\end{theorem}

{\bf Proof}   Taking 
$(c c^{\pi}) \in \pi$ , $(c c^{\pi})^{\mu}$ 
equals to $(c^{\pi} c^{\pi \varrho})$ or $(c^{\varrho}
c^{\pi \pi}) = (c^{\varrho} c)$ both cases giving $(c c^{\pi})^{\mu} \in
\varrho$.

By $\mu (\pi, \varrho)$  we denote arbitrary knot of some map belonging to 
$K_{\varrho}$  . Previous theorem and \ref{fix3} gives what follows .

\begin{theorem}
 $K_{\varrho} = \mu (\pi, \varrho) \times K_{\pi}$.  
\end{theorem}
This gives right for what follows.
\begin{corollary}
Every map can be expressed as a knot of this map multiplied by some 
selfconjugate map.
\end{corollary}

\section{Isomorphism}

Let us notice that  $A \in K_{\pi}$ if and only if $\pi^A = \pi$, because $\pi^{A^{-1}} =
\varrho_A$. Two maps $(P, P \cdot \pi_1)$ and $(Q, \cdot \pi_2)$ are isomorphic if and 
only if there exist such permutation $A$ that $\pi_2^A= \pi_1$ and 
$Q^A = P$. We write $P  \simeq Q$, saying that $P$ is {\em isomorphic} to 
$Q$ because they are conjugate with respect to  $A$.
\begin{theorem}
 Let $(\pi_1, P) \simeq (\pi_2,Q)$, i.e. they are conjugate
with respect to $A$. Both maps belong to one closed against
multiplication class of maps iff $A \in K_{\pi}$ with $\pi = \pi_1$.
\end{theorem}
{\bf Proof}   If $A \in K_{\pi_1}$ then $\pi_1^A$ is equal to $\pi_1$ ,but
because of conjugacy of $\pi_1$ and $\pi_2$ against $A$ also equal to
$\pi_2$. So $\pi_1 = \pi_2$ and $P$ and $Q$ are maps in one class with
$\pi = \pi_1 = \pi_2$. Conversely, if $\pi_1 = \pi_2$, then $\pi_1^A = \pi_2 = \pi_1$ and 
$A \in K_{\pi}$.

\section{Acknowledgements}
The author would like to thank Professor Ne\v{s}etr\'{\i}l and Dr. Kratochv\'{\i}l for 
giving the possibility to work at  the Caroline University in Prague in  spring 1994, during which
 time this work  mainly was done. The author would like to thank also Dr. 
U. Raitums for special care to get a grant for graph theoretical 
investigations. Many thanks also to Dr. Kratochv\'{\i}l, Prof. Ne\v{s}etr\'{\i}l and Dr. 
\c{K}ikusts for helpful discussions.

\end{document}